\documentclass{llncs}

\usepackage{amssymb}
\renewcommand{\pmod}[1]{\ (\mbox{mod }#1)}

\begin{document}
\bibliographystyle{splncs03}

\title{Sieving for pseudosquares and pseudocubes in parallel 
  using doubly-focused enumeration and wheel datastructures%
  \thanks{Supported by a grant from the Holcomb Awards Committe, and computing
   resources provided by the Frank Levinson Supercomputing Center
   at Butler University.}}

\author{Jonathan P.~Sorenson}%\inst{1}}

\institute{Butler University, Indianapolis IN 46208, USA, \\
\email{sorenson@butler.edu}, \\
\texttt{http://www.butler.edu/\homedir sorenson}
}

\maketitle

\begin{abstract}
  We extend the known tables of
  pseudosquares and pseudocubes,
  discuss the implications of these new data on the conjectured 
  distribution of pseudosquares and pseudocubes,
  and present the details of the algorithm used to do this work.
  Our algorithm is based on the space-saving wheel data structure
  combined with doubly-focused enumeration, run in parallel on a cluster
  supercomputer.
\end{abstract}
 
\section{Introduction}

It is well-known that testing for primality can be
  done in polynomial time \cite{AKS04,Bernstein07}.
However, the fastest known deterministic algorithms
  are conjectured to be the
  pseudosquares prime test of Lukes, Patterson, and Williams \cite{LPW96},
  and its generalization, 
  the pseudocube prime test of Berrizbeitia, M\"uller, and Willimas \cite{BMW04},
  both of which run in roughly cubic time,
  if a sufficiently large pseudosquare or pseudocube is available.
In particular, 
  the pseudosquares prime test is very useful in the context of finding
  all primes in an interval \cite{Sorenson06}, where sieving can be used
  in place of trial division.
This, then, motivates the search for larger and larger peudosquares
  and pseudocubes, and attempts to predict their distribution.
See, for example, Wooding and Williams \cite{WW2006} and also 
  \cite{PS2009,Williams,Schinzel97,Bernstein04,SW90}.

In this paper, we present extensions to the known tables
  of pseudosquares and pseudocubes in \S\ref{sec:results}.
We discuss the implications of this new data on the conjectured
  distribution of pseudosquares and pseudocubes in \S\ref{sec:distr},
  and give a minor refinement of the current conjectures.
Then we
  describe our parallel algorithm, based on 
  Bernstein's doubly-focused enumeration \cite{Bernstein04}, 
  which is used in a way similar, but not identical to the work of
  Wooding and Williams \cite{WW2006}, 
  combined with the space-saving wheel data structure 
  presented in \cite[\S4.1]{Sorenson06}.
We then suggest ideas for future work in \S\ref{sec:future}.

\section{Computational Results\label{sec:results}}

\newcommand{\pss}{L_{p,2}}
\newcommand{\psc}{L_{p,3}}
\newcommand{\pssn}{L_{p_n,2}}
\newcommand{\pscn}{L_{q_n,3}}

Let $(x/y)$ denote the Legendre symbol \cite{HW}.
For an odd prime $p$, let $\pss$, the \textit{pseudosquare} for $p$, be the
smallest positive integer such that
\begin{enumerate}
  \item $\pss \equiv 1 \pmod 8$,
  \item $(\pss / q)=1$ for every odd prime $q\le p$, and
  \item $\pss$ is not a perfect square.
\end{enumerate}
In other words, $\pss$ is a square modulo all primes up to $p$,
  but is not a square.
We found the following new pseudosquares:
\begin{quote}
\begin{tabular}{|ccc|crc|}
\hline
&$p$&& & $\pss$& \\ \hline\hline
&367&&& 36553 34429 47705 74600 46489& \\
&373&&& 42350 25223 08059 75035 19329& \\
&379&&& $> 10^{25}$& \\
\hline
\end{tabular}
\end{quote}
The two pseudosquares listed were found in 2008 in a computation that
went up to $5\times 10^{24}$, taking roughly 3 months wall time.
The final computation leading to the lower bound of $10^{25}$ ran for 
about 6 months, in two 3-month pieces, the second of which finished 
on January 1st, 2010.

Wooding and Williams \cite{WW2006} had found a lower bound of
  $L_{367,2} > 120120 \times 2^{64} \approx 2.216 \times 10^{24}$.
(Note: a complete table of pseudosquares, 
  current as of this writing, is available at
  \texttt{http://cr.yp.to/focus.html} care of Dan Bernstein).

Similarly, for an odd prime $p$, 
let $\psc$, the \textit{pseudocube} for $p$, be the
smallest positive integer such that
\begin{enumerate}
  \item $\psc \equiv \pm 1 \pmod 9$,
  \item $\psc^{(q-1)/3} \equiv1 \pmod q$ 
      for every prime $q\le p$, $q\equiv 1\pmod 3$, 
  \item $\gcd(\psc,q)=1$ for every prime $q\le p$, and
  \item $\psc$ is not a perfect cube.
\end{enumerate}
We found the following new pseudocubes (only listed for $p\equiv 1\pmod 3$):
\begin{quote}
\begin{tabular}{|ccc|crc|}
\hline
&$p$ &&& $\psc$& \\ \hline\hline
&499&&   &       601 25695 21674 16551 89317& \\
&523,541&&&     1166 14853 91487 02789 15947& \\
&547&    &&    41391 50561 50994 78852 27899& \\
&571,577&&&  1 62485 73199 87995 69143 39717& \\
&601,607&&&  2 41913 74719 36148 42758 90677& \\
&613    &&& 67 44415 80981 24912 90374 06633& \\
&619&&& $> 10^{27}$& \\
\hline
\end{tabular}
\end{quote}
These pseudocubes were found in about 6 months of total wall time in 2009.
Wooding and Williams \cite{WW2006} had found a lower bound of
  $L_{499,3} > 1.45152 \times 10^{22}$.
For a complete list of known pseudocubes, see
  \cite{WW2006,BMW04,SW90}.

\section{The Distribution of Pseudosquares and Pseudocubes\label{sec:distr}}

Let $p_i$ denote the $i$th prime, and $q_i$ denote the $i$th prime such that
  $q_i \equiv 1 \pmod 3$.
In \cite{LPW96} it was conjectured that, for a constant $c_2>0$, we have
\begin{equation}\label{eq:pssdist}
   \pssn \approx c_2 2^n \log p_n.
\end{equation}
Using similar methods, in \cite{BMW04} it was conjectured that, 
  for a constant $c_3>0$, we have
\begin{equation}\label{eq:pscdist}
   \pscn \approx c_3 3^n (\log q_n)^2.
\end{equation}

\begin{table}
\caption{Values of $c_2(n)$ based on known pseudosquares.\label{tbl:pss}}
\begin{quote}
\begin{tabular}{rc|crc|crc|cr}
  $n$ &&& $p_n$ &&& $\pssn$ &&& $c_2(n)$ \\ \hline
2 &&& 3 &&& 73 &&& 16.61 \\
3 &&& 5 &&& 241 &&& 18.72 \\
4 &&& 7 &&& 1009 &&& 32.41 \\
5 &&& 11 &&& 2641 &&& 34.42 \\
6 &&& 13 &&& 8089 &&& 49.28 \\
7 &&& 17 &&& 18001 &&& 49.64 \\
8 &&& 19 &&& 53881 &&& 71.48 \\
9 &&& 23 &&& 87481 &&& 54.49 \\
10 &&& 29 &&& 117049 &&& 33.95 \\
11 &&& 31 &&& 515761 &&& 73.34 \\
12 &&& 37 &&& 1083289 &&& 73.24 \\
13 &&& 41 &&& 3206641 &&& 105.41 \\
14 &&& 43 &&& 3818929 &&& 61.97 \\
15 &&& 47 &&& 9257329 &&& 73.38 \\
16 &&& 53 &&& 22000801 &&& 84.55 \\
17 &&& 59 &&& 48473881 &&& 90.70 \\
18 &&& 61 &&& 48473881 &&& 44.98 \\
19 &&& 67 &&& 175244281 &&& 79.49 \\
20 &&& 71 &&& 427733329 &&& 95.70 \\
21 &&& 73 &&& 427733329 &&& 47.54 \\
22 &&& 79 &&& 898716289 &&& 49.04 \\
23 &&& 83 &&& 2805544681 &&& 75.69 \\
24 &&& 89 &&& 2805544681 &&& 37.25 \\
25 &&& 97 &&& 2805544681 &&& 18.28 \\
26 &&& 101 &&& 10310263441 &&& 33.29 \\
27 &&& 103 &&& 23616331489 &&& 37.96 \\
28 &&& 107 &&& 85157610409 &&& 67.89 \\
29 &&& 109 &&& 85157610409 &&& 33.81 \\
30 &&& 113 &&& 196265095009 &&& 38.67 \\
31 &&& 127 &&& 196265095009 &&& 18.87 \\
32 &&& 131 &&& 2871842842801 &&& 137.15 \\
33 &&& 137 &&& 2871842842801 &&& 67.95 \\
34 &&& 139 &&& 2871842842801 &&& 33.88 \\
35 &&& 149 &&& 26250887023729 &&& 152.68 \\
36 &&& 151 &&& 26250887023729 &&& 76.14 \\
37 &&& 157 &&& 112434732901969 &&& 161.79 \\
38 &&& 163 &&& 112434732901969 &&& 80.30 \\
\end{tabular} $\quad$
\begin{tabular}{rc|crc|crc|cr}
  $n$ &&& $p_n$ &&& $\pssn$ &&& $c_2(n)$ \\ \hline
39 &&& 167 &&& 112434732901969 &&& 39.96 \\
40 &&& 173 &&& 178936222537081 &&& 31.58 \\
41 &&& 179 &&& 178936222537081 &&& 15.69 \\
42 &&& 181 &&& 696161110209049 &&& 30.45 \\
43 &&& 191 &&& 696161110209049 &&& 15.07 \\
44 &&& 193 &&& 2854909648103881 &&& 30.84 \\
45 &&& 197 &&& 6450045516630769 &&& 34.70 \\
46 &&& 199 &&& 6450045516630769 &&& 17.32 \\
47 &&& 211 &&& 11641399247947921 &&& 15.46 \\
48 &&& 223 &&& 11641399247947921 &&& 7.65 \\
49 &&& 227 &&& 190621428905186449 &&& 62.42 \\
50 &&& 229 &&& 196640148121928601 &&& 32.14 \\
51 &&& 233 &&& 712624335095093521 &&& 58.06 \\
52 &&& 239 &&& 1773855791877850321 &&& 71.92 \\
53 &&& 241 &&& 2327687064124474441 &&& 47.12 \\
54 &&& 251 &&& 6384991873059836689 &&& 64.15 \\
55 &&& 257 &&& 8019204661305419761 &&& 40.11 \\
56 &&& 263 &&& 10198100582046287689 &&& 25.40 \\
57 &&& 269 &&& 10198100582046287689 &&& 12.65 \\
58 &&& 271 &&& 10198100582046287689 &&& 6.32 \\
59 &&& 277 &&& 69848288320900186969 &&& 21.54 \\
60 &&& 281 &&& 208936365799044975961 &&& 32.14 \\
61 &&& 283 &&& 533552663339828203681 &&& 40.99 \\
62 &&& 293 &&& 936664079266714697089 &&& 35.76 \\
63 &&& 307 &&& 936664079266714697089 &&& 17.73 \\
64 &&& 311 &&& 2142202860370269916129 &&& 20.23 \\
65 &&& 313 &&& 2142202860370269916129 &&& 10.10 \\
66 &&& 317 &&& 2142202860370269916129 &&& 5.04 \\
67 &&& 331 &&& 13649154491558298803281 &&& 15.94 \\
68 &&& 337 &&& 34594858801670127778801 &&& 20.14 \\
69 &&& 347 &&& 99492945930479213334049 &&& 28.81 \\
70 &&& 349 &&& 99492945930479213334049 &&& 14.39 \\
71 &&& 353 &&& 295363187400900310880401 &&& 21.32 \\
72 &&& 359 &&& 295363187400900310880401 &&& 10.63 \\
73 &&& 367 &&& 3655334429477057460046489 &&& 65.54 \\
74 &&& 373 &&& 4235025223080597503519329 &&& 37.86 \\
\  \\%& && &&& &&&
\end{tabular}
\end{quote}
\end{table}
\begin{table}
\caption{Values of $c_3(n)$ based on known pseudocubes.\label{tbl:psc}}
\begin{quote}
\begin{tabular}{rc|crc|crc|cl}
  $n$ &&& $q_n$ &&& $\pscn$ &&& $c_3(n)$ \\ \hline
10 &&& 79 &&& 7235857 &&& 6.42 \\
11 &&& 97 &&& 8721539 &&& 2.35 \\
12 &&& 103 &&& 8721539 &&& 0.764 \\
13 &&& 109 &&& 91246121 &&& 2.6 \\
14 &&& 127 &&& 91246121 &&& 0.813 \\
15 &&& 139 &&& 98018803 &&& 0.281 \\
16 &&& 151 &&& 1612383137 &&& 1.49 \\
17 &&& 157 &&& 1612383137 &&& 0.488 \\
18 &&& 163 &&& 7991083927 &&& 0.795 \\
19 &&& 181 &&& 7991083927 &&& 0.254 \\
20 &&& 193 &&& 7991083927 &&& 0.0827 \\
21 &&& 199 &&& 20365764119 &&& 0.0695 \\
22 &&& 211 &&& 2515598768717 &&& 2.8 \\
23 &&& 223 &&& 6440555721601 &&& 2.34 \\
24 &&& 229 &&& 29135874901141 &&& 3.49 \\
25 &&& 241 &&& 29135874901141 &&& 1.14 \\
26 &&& 271 &&& 29135874901141 &&& 0.365 \\
27 &&& 277 &&& 406540676672677 &&& 1.69 \\
28 &&& 283 &&& 406540676672677 &&& 0.558 \\
29 &&& 307 &&& 406540676672677 &&& 0.181 \\
30 &&& 313 &&& 406540676672677 &&& 0.0598 \\
31 &&& 331 &&& 75017625272879381 &&& 3.61 \\
\end{tabular} $\quad$
\begin{tabular}{rc|crc|crc|cl}
  $n$ &&& $q_n$ &&& $\pscn$ &&& $c_3(n)$ \\ \hline
32 &&& 337 &&& 75017625272879381 &&& 1.2 \\
33 &&& 349 &&& 75017625272879381 &&& 0.394 \\
34 &&& 367 &&& 996438651365898469 &&& 1.71 \\
35 &&& 373 &&& 2152984914389968651 &&& 1.23 \\
36 &&& 379 &&& 12403284862819956587 &&& 2.34 \\
37 &&& 397 &&& 37605274105479228611 &&& 2.33 \\
38 &&& 409 &&& 37605274105479228611 &&& 0.77 \\
39 &&& 421 &&& 37605274105479228611 &&& 0.254 \\
40 &&& 433 &&& 205830039006337114403 &&& 0.459 \\
41 &&& 439 &&& 1845193818928603436441 &&& 1.37 \\
42 &&& 457 &&& 7854338425385225902393 &&& 1.91 \\
43 &&& 463 &&& 12904554928068268848739 &&& 1.04 \\
44 &&& 487 &&& 13384809548521227517303 &&& 0.355 \\
45 &&& 499 &&& 60125695216741655189317 &&& 0.527 \\
46 &&& 523 &&& 116614853914870278915947 &&& 0.336 \\
47 &&& 541 &&& 116614853914870278915947 &&& 0.111 \\
48 &&& 547 &&& 4139150561509947885227899 &&& 1.31 \\
49 &&& 571 &&& 16248573199879956914339717 &&& 1.69 \\
50 &&& 577 &&& 16248573199879956914339717 &&& 0.56 \\
51 &&& 601 &&& 24191374719361484275890677 &&& 0.274 \\
52 &&& 607 &&& 24191374719361484275890677 &&& 0.0912 \\
53 &&& 613 &&& 674441580981249129037406633 &&& 0.845 
\end{tabular}
\end{quote}
\end{table}

In a desire to test the accuracy of these conjectures,
  for integers $n>0$ let us define
\begin{eqnarray}
  c_2(n) &:=& \frac{\pssn}{ 2^n \log p_n}, \\
  c_3(n) &:=& \frac{\pscn}{  3^n (\log q_n)^2}.
\end{eqnarray}
We calculated $c_2(n)$ and $c_3(n)$ from known pseudosquares and
  pseudocubes.
We present these computations in Table \ref{tbl:pss}, for pseudosquares,
  and in Table \ref{tbl:psc}, for pseudocubes, below.

From Table \ref{tbl:pss}, we readily see that
  $c_2(n)$ appears to be bounded between roughly
  $5$ and $162$, with an average value near $45$.
There is no clear trend toward zero or infinity.
Due to the common occurence of values of $n$ where
  $L_{p_n,2}=L_{p_{n+1},2}$ (for example, $n=56$), 
  it should also be clear $c_2(n)$ does not have a limit.

Similarly for the pseudocubes, in Table \ref{tbl:psc} we see
  that $0.05 < c_3(n) < 6.5$ for $10\le n \le 53$,
  with an average value of roughly $1.22$.
And again, there is no clear trend toward zero or infinity, nor can there
  be a limit for $c_3(n)$.

This leads us to the following refinements, if you will, of the
  conjectures (\ref{eq:pssdist}),(\ref{eq:pscdist}) above.
\paragraph{Conjecture.}
  For the pseudosquares, we conjecture that
 \begin{eqnarray}
   \liminf_{n\rightarrow\infty} \frac{\pssn}{2^n\log p_n} &>& 0, \\
   \limsup_{n\rightarrow\infty} \frac{\pssn}{2^n\log p_n} &<& \infty.
  \end{eqnarray}
Similarly, for the pseudocubes, we conjecture that
 \begin{eqnarray}
   \liminf_{n\rightarrow\infty} \frac{\pscn}{3^n(\log q_n)^2} &>& 0, \\
   \limsup_{n\rightarrow\infty} \frac{\pscn}{3^n(\log q_n)^2} &<& \infty.
  \end{eqnarray}

Our data also has implications on the relative efficiently of
  primality testing.
In particular, several researchers have pointed out that if
  conjectures (\ref{eq:pssdist}),(\ref{eq:pscdist}) are true,
  then the running time of the pseudocube prime test,
  which depends on the value of $\pscn^{2/3}$,
  should eventually outperform the pseudosquare prime test,
  whose running time depends on $\pssn$.
In particular, one infers from 
  conjectures (\ref{eq:pssdist}) and (\ref{eq:pscdist}) that
\begin{equation}\label{eq:faster}
 \frac{\pscn^{2/3}}{\pssn}
  \quad \gg \quad \left( \frac{3^{2/3}}{2} \right)^n \quad > \quad 1 
\end{equation}
for sufficiently large $n$ (see \cite[\S9.1]{WW2006}).
This inference follows from our refined conjectures as well.

We have our first specific value of $n$
  to support (\ref{eq:faster}),
  namely with $n=48$,
  where $ {\pscn^{2/3}} \approx 2.214 \cdot \pssn$.
However, given that $c_2(n)$ averages about 45,
  and $c_3(n)$ averages just over 1.2, 
  we would reasonably expect (\ref{eq:faster}) to largely be true only for $n$ 
  larger than about $75$, under the assumption these averages are maintained.
To test this,
  more pseudosquares and, in particular, more pseudocubes are needed.

\section{Algorithm Details\label{sec:alg}}

We begin with a review of doubly-focused enumeration,
  explain how we employ parallelism,
  and how the space-saving wheel datastructure is utilized.
We also discuss the details of our implementation, including the
  hardware platform and software used.

\subsection{Doubly-Focused Enumeration}

The main idea is that every integer $x$, with $0\le x\le H$,
can be written in the form
\begin{equation}
 x=t_pM_n - t_nM_p
\end{equation}
where
\begin{equation}
 \gcd(M_p,M_n)=1, \quad 0\le t_p \le \frac{H+M_nM_p}{M_n}, \quad
    \mbox{and} \quad 0\le t_n < M_n.
\end{equation}
(See \cite{Bernstein04} or \cite[Lemma 1]{WW2006}.)
This is an explicit version of the Chinese Remainder Theorem.

To find pseudosquares, we set $M_n$ and $M_p$ to be products of small odd
primes and $8$, choose $t_p$ to be square modulo $M_p$, 
and $-t_n$ to be square modulo $M_n$.
To be precise, in our implementation we set
\begin{eqnarray*}
  M_p&=&  7\cdot 11\cdot 13\cdot 17\cdot 19\cdot 23\cdot 29\cdot 31
          \cdot 37\cdot 41\cdot 43\cdot 53\cdot 89 \\
    &=&2057\, 04617\, 33829\, 17717 \qquad \mbox{and} \\
  M_n&=&  8\cdot 3\cdot 5\cdot 47\cdot 59\cdot 61\cdot 67\cdot 71\cdot 73
           \cdot 79\cdot 83\cdot 97 \\
    &=&4483\, 25952\, 77215\, 26840.
\end{eqnarray*}
Note that both $M_p,M_n < 2^{64}$, allowing us to work in 64-bit 
  machine arithmetic. 
%Including one more prime (101) would force one of $M_p,M_n$ to exceed $2^{64}$.

To find pseudocubes, the same idea applies, only note that if $-t_n$ is a
  cube modulo $M_n$, so is $t_n$.  We used only $2,9$ and primes congruent to
  $1 \pmod 3$ for better filter rates:
\begin{eqnarray*}
  M_p&=& 2\cdot 7\cdot 13\cdot 31\cdot 43\cdot 73\cdot 79\cdot 127\cdot 139
          \cdot 157\cdot 181 \\
    &=&701\, 85635\, 61110\, 39402 
 \qquad \mbox{and} \\
  M_n&=& 9\cdot 19\cdot 37\cdot 61\cdot 67\cdot 97\cdot 103\cdot 109
          \cdot 151\cdot 163 \\
    &=&693\, 11050\, 43291\, 92503
\end{eqnarray*}

\subsection{Parallelism and Main Loop}

Each processor core was assigned an interval of $t_p$ values to process
  by giving it values of $H^{-}$ and $H^{+}$.

For finding pseudosquares, $H^+-H^- \approx M_n \cdot 4.76 \times 10^{11}$.
For finding pseudocubes, $H^+-H^-  \approx  M_n \cdot 4.99 \times 10^{12}$.

Parallelism was achieved by having different processors working on
  different intervals simultaneously.
Once all processors had finished their current intervals,
  the work was saved to disk
  (allowing restarts as needed) and new intervals were assigned.

To process an interval, each processor core did the following:
\begin{enumerate}
\item
  Using the wheel datastructure, generate all square or cube values of
    $t_p$ with $H^- \le t_pM_n \le H^+$, and store these in an array 
    \texttt{A[]}.
\item
  The wheel datastructure does not generate the $t_p$ values in order, so
    sort \texttt{A[]} in memory using quicksort.
  Note that $H^-$ and $H^+$ are chosen close enough together so that
    this array held no more than 40 million integers, using
    at most 320 megabytes of RAM per processor core.

\item
  Using the first and last entries in \texttt{A[]}, compute a range of
    valid $t_n$ values to process,
  and then use a wheel datastructure to generate all $t_n$ values in that range
    such that $-t_n$ is square modulo $M_n$ for pseudosquares,
    or $t_n$ is a cube modulo $M_n$ for pseudocubes.

  We use an outer loop over $t_n$ values in the order enumerated by the
    wheel data structure for $M_n$, and an inner loop over
    consecutive $t_p$ values drawn from \texttt{A[]}.

\item
  For each $t_n$ generated, we normalize sieve tables for the next
    4 primes ($101,103,107,109$ for pseudosquares, and
    $193,199,211,223$ for pseudocubes) to allow for constant-time
    table lookup to see if an $x$-value (see below) is a square/cube
    modulo these primes, indexed by $t_p$ value.

    The number of primes to use for this depends on how many $t_p$ values
      will be processed for each $t_n$ -- in our case, it was several
      hundred on average, so this step improves performance.
    If it were fewer, say 50, then normalizing the sieve tables would require more
      work than is saved by having constant-time lookup.

\item
  For each $t_n$ generated, using binary search on \texttt{A[]} to find
    all the $t_p$ values it can match with, generate an
    $ x=t_pM_n - t_nM_p$ within our global search range.
  (For example, in our last run for pseudosquares, we searched for $x$
    values between $7.5 \times 10^{24}$ and $10^{25}$.)

   Note: at this point we do not actually compute the value of $x$.

\item
  Lookup each $t_p$ value in the normalized tables mentioned above.
  If it fails any of the 4 sieve tests, move on to the next $t_p$ value.
  For pseudosquares, a $t_p$ values passes these tests with probability
   roughly $(1/2)^4 = 1/16$, and for pseudocubes, roughly $(1/3)^4 = 1/81$.

  Note that this step is the running time bottleneck of the algorithm.

\item
  The next batch of primes $q$ have precomputed sieve tables that are not
   normalized, but we precompute $M_p$ and $M_n$ modulo each $q$
   so the we can compute $x \bmod q$ without exceeding 64-bit arithmetic.
  Continue only if our $t_p$ value passes all these sieve tests as well.
  The expected number of primes $q$ used in this step is constant.

\item
  Finally, compute $x$ using 128-bit hardware arithmetic, and see if it
    is a perfect square or perfect cube.
  If it passes this test, append $x$ to the output file 
    for this processor core.
  
\end{enumerate}
We had two wheel datastructures, one each for $M_p$ and $M_n$.
For details on how this datastructure works, see \cite{Sorenson06}.
We leave the details for how to modify the datastructure to handle
  cubes in place of squares to the reader.

\subsection{Implementation Details}

To compute the tables presented in \S\ref{sec:results}, 
  we used Butler University's cluster supercomputer,
  \textit{BigDawg}, which has 24 compute nodes, each of which has
  four AMD Opteron 8354 quad-core CPUs at 2.2GHz with 512KB cache, 
  for a total of 384 compute cores.
As might be expected, we did not have sole access to this machine for over
  a year, so the code was designed, and ran, using anywhere from 10 to 24 nodes,
  or from 160 to 384 cores,
  depending on the needs of other users.
This flexibility is one advantage of our parallelization method -- by 
  $t_p$ intervals.
In \cite{WW2006}, they parallelized over residue classes, which restricts
  the CPU count to a fixed number (180 in their case).

\textit{BigDawg} runs a Linux kernel on its head node and compute nodes,
  and the code was written in C++ using the gnu compiler (version 4.1.2) 
  with MPI.
It has both 10GB ethernet and
  Infiniband interconnect, but inter-processor communication was
  not a bottleneck for our programs.

We tested our code by first finding known pseudosquares 
  (all but the highest few) and known pseudocubes, in the process verifying
  previous results.

\section{Future Work\label{sec:future}}

We plan to port our code to work with 8 NVidia GPUs recently added to
Butler's supercomputer, giving it roughly 2-3 times the raw computing power.
This will require a major restructuring of the code, and the removal of
recursion in the wheel datastructure.

%\bibliography{all}

\end{document}